\documentclass[12pt]{article}%
\usepackage{amsmath}
\usepackage{amsfonts}
\usepackage{amssymb}
\usepackage{graphicx}%
\providecommand{\U}[1]{\protect\rule{.1in}{.1in}}
\newtheorem{theorem}{Theorem}[section]

\newtheorem{corollary}[theorem]{Corollary}

\newtheorem{definition}{Definition}[section]

\newtheorem{observation}[theorem]{Observation}

\newtheorem{proposition}[theorem]{Proposition}

\providecommand{\boksie}{\ensuremath{\mathbin{\raisebox{0.3mm}{$\scriptstyle\square$}}}}
\begin{document}

\title{\textbf{Boundary Independent Broadcasts}\\\textbf{in Graphs}}
\author{C.M. Mynhardt\thanks{Supported by the Natural Sciences and Engineering
Research Council of Canada.} and L. Neilson\thanks{Department of Adult Basic
Education, Vancouver Island University, 900 Fifth Street, Nanaimo, BC;
\textsc{Canada} V9R 5S5}\\Department of Mathematics and Statistics\\University of Victoria, Victoria, BC, \textsc{Canada}\\{\small kieka@uvic.ca, linda.neilson@viu.ca\bigskip}\\{\small This paper is dedicated to our friend and colleague, }\\{\small Gary MacGillivray, on the occasion of his 60}$^{\mathrm{th}}%
${\small birthday. }\\{\small Thanks, Gary, for being there for your colleagues and students!}}
\date{}
\maketitle

\begin{abstract}
A broadcast on a nontrivial connected graph $G=(V,E)$ is a function
$f:V\rightarrow\{0,1,\dots,\operatorname{diam}(G)\}$ such that $f(v)\leq e(v)$
(the eccentricity of $v$) for all $v\in V$. The weight of $f$ is $\sigma(f)=%
{\textstyle\sum_{v\in V}}
f(v)$. A vertex $u$ hears $f$ from $v$ if $f(v)>0$ and $d(u,v)\leq f(v)$. A
broadcast $f$ is independent, or hearing independent, if no vertex $u$ with
$f(u)>0$ hears $f$ from any other vertex $v$. We define a different type of
independent broadcast, namely a boundary independent broadcast, as a broadcast
$f$ such that, if a vertex $w$ hears $f$ from vertices $v_{1},...,v_{k}%
,\ k\geq2$, then $d(w,v_{i})=f(v_{i})$ for each $i$. The maximum weights of a
hearing independent broadcast and a boundary independent broadcast are the
hearing independence broadcast\emph{ }number $\alpha_{h}(G)$ and the boundary
independence broadcast number\emph{ }$\alpha_{\operatorname{bn}}(G)$, respectively.

We prove that $\alpha_{\operatorname{bn}}(G)=\alpha(G)$ (the independence
number) for any $2$-connected bipartite graph $G$ and that $\alpha
_{\operatorname{bn}}(G)\leq n-1$ for all graphs $G$ of order $n$,
characterizing graphs for which equality holds. We compare $\alpha
_{\operatorname{bn}}$ and $\alpha_{h}$ and prove that although the difference
$\alpha_{h}-\alpha_{\operatorname{bn}}$ can be arbitrary, the ratio is
bounded, namely $\alpha_{h}/\alpha_{\operatorname{bn}}<2$, which is
asymptotically best possible. We deduce that $\alpha_{h}(G)\leq2n-5$ for all
connected graphs $G\neq P_{n}$ of order $n$, which improves an existing upper
bound for $\alpha_{h}(G)$ when $\alpha(G)\geq n/2$.

\end{abstract}

\noindent\textbf{Keywords:\hspace{0.1in}}broadcast domination; broadcast
independence, hearing independence; boundary independence

\noindent\textbf{AMS Subject Classification Number 2010:\hspace{0.1in}}05C69

\section{Introduction}

In a search for the best way to generalize the concept of independent sets in
graphs to independent broadcasts, there are several ways to look at an
independent set $X$ of a graph $G$. One way is from the point of view of the
vertices in $X$: no two vertices are adjacent -- the usual definition. Another
way is from the point of view of the edges of $G$: no edge is incident with
(or covered by) more than one vertex in $X$. Using the latter approach we
define boundary independent broadcasts as an alternative to independent
broadcasts as defined by Erwin \cite{Ethesis}, which we refer to here as
hearing independent broadcasts. Among other results we show that the boundary
independent broadcast number $\alpha_{\operatorname{bn}}$ of any graph lies
between its independence number and its hearing independent broadcast number
$\alpha_{h}$. We prove a tight upper bound for $\alpha_{\operatorname{bn}}$
which leads to a new tight upper bound for $\alpha_{h}$.

\subsection{Broadcast definitions}

\label{SecBC}For undefined concepts we refer the reader to \cite{CLZ}. The
study of broadcast domination was initiated by Erwin in his doctoral
dissertation \cite{Ethesis}. A \emph{broadcast} on a nontrivial connected
graph $G=(V,E)$ is a function $f:V\rightarrow\{0,1,\dots,\operatorname{diam}%
(G)\}$ such that $f(v)\leq e(v)$ (the eccentricity of $v$) for all $v\in V$.
When $G$ is disconnected, we define a broadcast on $G$ as the union of
broadcasts on its components. Define $V_{f}^{+}=\{v\in V:f(v)>0\}$ and
partition $V_{f}^{+}$ into the two sets $V_{f}^{1}=\{v\in V:f(v)=1\}$ and
$V_{f}^{++}=V_{f}^{+}-V_{f}^{1}$. A vertex in $V_{f}^{+}$ is called a
\emph{broadcasting vertex}. A vertex $u$ \emph{hears} $f$ from $v\in V_{f}%
^{+}$, and $v$ $f$-\emph{dominates} $u$, if the distance $d(u,v)\leq f(v)$. If
$d(u,v)<f(v)$, we also say that say that $v$ \emph{overdominates }$u$. Denote
the set of all vertices that do not hear $f$ by $U_{f}$. A broadcast $f$ is
\emph{dominating} if $U_{f}=\varnothing$. The \emph{weight} of $f$ is
$\sigma(f)=\sum_{v\in V}f(v)$, and the \emph{broadcast number} of $G$ is
\[
\gamma_{b}(G)=\min\left\{  \sigma(f):f\text{ is a dominating broadcast of
}G\right\}  .
\]
When $f$ and $g$ are broadcasts on $G$ such that $g(v)\leq f(v)$ for each
$v\in V$, we write $g\leq f$. When in addition $g(v)<f(v)$ for at least one
$v\in V$, we write $g<f$. A dominating broadcast $f$ on $G$ is a \emph{minimal
dominating broadcast} if no broadcast $g<f$ is dominating. The \emph{upper
broadcast number }of $G$ is
\[
\Gamma_{b}(G)=\max\left\{  \sigma(f):f\text{ is a minimal dominating broadcast
of }G\right\}  ,
\]
and a dominating broadcast $f$ of $G$ such that $\sigma(f)=\Gamma_{b}(G)$ is
called a $\Gamma_{b}$-\emph{broadcast}. First defined by Erwin \cite{Ethesis},
the upper broadcast number was also studied by Ahmadi, Fricke, Schroeder,
Hedetniemi and Laskar \cite{Ahmadi}, Bouchemakh and Fergani \cite{BF}, Dunbar,
Erwin, Haynes, Hedetniemi and Hedetniemi \cite{Dunbar}, Gemmrich and Mynhardt
\cite{GM} and Mynhardt and Roux \cite{MR}.

If $f$ is a (minimal) dominating broadcast such that $V_{f}^{+}=V_{f}^{1}$,
then $f$ is the characteristic function of a (minimal) dominating set. Hence,
denoting the cardinalities of a minimum dominating set and a maximum minimal
dominating set by $\gamma(G)$ and $\Gamma(G)$ (the \emph{lower} and
\emph{upper domination numbers} of $G$), respectively, we see that $\gamma
_{b}(G)\leq\gamma(G)$ and $\Gamma(G)\leq\Gamma_{b}(G)$ for any graph $G$.

We denote the independence number of $G$ by $\alpha(G)$ and the minimum
cardinality of a maximal independent set (the \emph{independent domination
number} of $G$) by\emph{ }$i(G)$. To generalize the concept of independent
sets, Erwin \cite{Ethesis} defined a broadcast $f$ to be \emph{independent},
or, for our purposes, \emph{hearing independent}, if no vertex $u\in V_{f}%
^{+}$ hears $f$ from any other vertex $v\in V_{f}^{+}$; that is, broadcasting
vertices only hear themselves. This version of broadcast independence was also
considered by, among others, Ahmane, Bouchemakh and Sopena \cite{ABS}, Bessy
and Rautenbach \cite{BR}, and Bouchemakh and Zemir \cite{Bouch}. We show below
that other definitions of broadcast independence, which also generalize
independent sets and lead to different independent broadcast numbers, are feasible.

\subsection{Neighbourhoods and boundaries}

Following \cite{MR}, for a broadcast $f$ on $G$ and $v\in V_{f}^{+}$, we
define the

\begin{itemize}
\item $f$-\emph{neighbourhood}$\ $of $v$ by $N_{f}(v)=\{u\in V:d(u,v)\leq
f(v)\}$,

\item $f$-\emph{boundary} of $v$ by $B_{f}(v)=\{u\in V:d(u,v)=f(v)\}$,

\item $f$-\emph{private neighbourhood}$\ $of $v$ by $\operatorname{PN}%
_{f}(v)=\{u\in N_{f}(v):u\notin N_{f}(w)$ for all$\ w\in V_{f}^{+}-\{v\}\}$,

\item $f$-\emph{private boundary}$\ $of $v$ by $\operatorname{PB}%
_{f}(v)=\{u\in N_{f}(v):u$ is not dominated by$\ (f-\{(v,f(v)\})\cup
\{(v,f(v)-1)\}$.
\end{itemize}

Note that if $u\in V_{f}^{1}$ and $u$ does not hear $f$ from any vertex $v\in
V_{f}^{+}-\{u\}$, then $u\in\operatorname{PB}_{f}(u)$, and if $u\in V_{f}%
^{++}$, then $\operatorname{PB}_{f}(u)=B_{f}(u)\cap\operatorname{PN}_{f}(u)$.
If $f$ is a broadcast such that every vertex $x$ that hears more than one
broadcasting vertex also satisfies $d(x,u)\geq f(u)$ for all $u\in V_{f}^{+}$,
then the \emph{broadcast only overlaps in boundaries}. On the other hand, if
$f$ is a dominating broadcast such that no vertex hears more than one
broadcasting vertex, then $f$ is an \emph{efficient dominating broadcast}.
When $xy\in E(G)$ and $x,y\in N_{f}(u)$ for some $u\in V_{f}^{+}$ such that at
least one of $x$ and $y$ does not belong to $B_{f}(u)$, we say that the edge
$xy$ is \emph{covered} in $f$ by $u$. When $xy$ is not covered by any $u\in
V_{f}^{+}$, we say that $xy$ is \emph{uncovered by~}$f$.

Erwin \cite{Ethesis} determined a necessary and sufficient condition for a
dominating broadcast to be minimal dominating. We restate it here in terms of
private boundaries.

\begin{proposition}
\label{PropMinimal}\emph{\cite{Ethesis}}\hspace{0.1in}A dominating broadcast
$f$ is a minimal dominating broadcast if and only if $\operatorname{PB}%
_{f}(v)\neq\varnothing$ for each $v\in V_{f}^{+}$.
\end{proposition}

Ahmadi et al.~\cite{Ahmadi} define \label{Def_ir}a broadcast $f$ to be
\emph{irredundant} if $\operatorname{PB}_{f}(v)\neq\varnothing$ for each $v\in
V_{f}^{+}$. An irredundant broadcast $f$ is \emph{maximal irredundant} if no
broadcast $g>f$ is irredundant. The \emph{lower} and \emph{upper broadcast
irredundant numbers} of $G$ are
\[
\operatorname{ir}_{b}(G)=\min\left\{  \sigma(f):f\text{ is a maximal
irredundant broadcast of }G\right\}
\]
and
\[
\operatorname{IR}_{b}(G)=\max\left\{  \sigma(f):f\text{ is an irredundant
broadcast of }G\right\}  ,
\]
respectively. Proposition \ref{PropMinimal} and the above definitions imply
the following two results.

\begin{corollary}
\label{Cor_ir-dom}\emph{\cite{Ahmadi}}\hspace{0.1in}$(i)\hspace{0.1in}$Any
minimal dominating broadcast is maximal\newline irredundant.

\begin{enumerate}
\item[$(ii)$] For any graph $G$,
\begin{equation}
\operatorname{ir}_{b}(G)\leq\gamma_{b}(G)\leq\gamma(G)\leq i(G)\leq
\alpha(G)\leq\Gamma(G)\leq\Gamma_{b}(G)\leq\operatorname{IR}_{b}(G).
\label{sequence1}%
\end{equation}

\end{enumerate}
\end{corollary}

\subsection{Independent broadcasts}

The characteristic function of an independent set has the following features,
which we generalize to obtain three different types of broadcast
independence:\smallskip

(a)$\hspace{0.1in}$\emph{boundary} or bn-independent type: broadcasts overlap
only in boundaries.

(b)\hspace{0.1in}\emph{hearing} or h-independent type \cite{Ethesis}:
broadcasting vertices hear only themselves.

(c)$\hspace{0.1in}$\emph{set} or s-independent type: broadcasting vertices
form an independent set.\smallskip

Broadcasts of type (c) were considered by Neilson \cite{LindaD} and found to
be not very interesting. We now consider broadcasts of type (a) and define
three new types of broadcast independence. Additional types can be found in
\cite{LindaD}. If a broadcast $f$ satisfies one of our definitions of
independence and there is no broadcast $g$ such that $g>f$ and $g$ also meets
our definition of independence, we say that $f$ is a \emph{maximal independent
broadcast} for this type of independence. Otherwise $f$ is not maximal
independent and can be \emph{extended} to a larger weight broadcast (for
example to $g$) which satisfies the given definition of
independence.\label{defs}

\begin{definition}
\label{bn-i}\emph{\cite{LindaD}\hspace{0.1in}A broadcast is} bn-independent
\emph{if it overlaps only in boundaries. The maximum (minimum) weight of a
(maximal) bn-independent broadcast on }$G$\emph{ is }$\alpha
_{\operatorname{bn}}(G)$ \emph{(}$i_{\operatorname{bn}}(G)$\emph{); such a
broadcast is called an }$\alpha_{\operatorname{bn}}$-broadcast\emph{
(}$i_{\operatorname{bn}}$-broadcast\emph{).}
\end{definition}

\begin{definition}
\label{bnr-i}\emph{\cite{LindaD}\hspace{0.1in}A broadcast is}
bnr-independent\emph{ if it is bn-independent and irredundant. The maximum
(minimum) weight of a (maximal) bnr-independent broadcast is} $\alpha
_{\operatorname{bnr}}(G)$ \emph{(}$i_{\operatorname{bnr}}(G)$\emph{); such a
broadcast is called an }$\alpha_{\operatorname{bnr}}$-broadcast\emph{
(}$i_{\operatorname{bnr}}$-broadcast\emph{).}
\end{definition}

\begin{definition}
\label{bnd-i}\emph{\cite{LindaD}\hspace{0.1in}A broadcast is}
bnd-independent\emph{ if it is minimal dominating and bn-independent. The
maximum (minimum) weight of a bnd-independent broadcast is} $\alpha
_{\operatorname{bnd}}(G)$ \emph{(}$i_{\operatorname{bnd}}(G)$\emph{); such a
broadcast is called an }$\alpha_{\operatorname{bnd}}$-broadcast\emph{
(}$i_{\operatorname{bnd}}$-broadcast\emph{).}
\end{definition}

\begin{definition}
\label{h-i}\emph{\cite{Ethesis}\hspace{0.1in}The maximum (minimum) weight of a
(maximal) h-independent broadcast is} $\alpha_{h}(G)$ \emph{(}$i_{h}%
(G)$\emph{); such a broadcast is called an }$\alpha_{h}$-broadcast\emph{
(}$i_{h}$-broadcast\emph{).}
\end{definition}

A bnd-independent broadcast, because it is minimal dominating, is maximal
irredundant (Corollary \ref{Cor_ir-dom}), and because it is irredundant and
dominating, it is minimal dominating (Proposition \ref{PropMinimal}). The
parameters $\alpha_{h}(G)$ and $\alpha_{\operatorname{bn}}(G)$ are also called
the \emph{hearing }or\emph{ h-independence broadcast number }and the
\emph{boundary} or\emph{ bn-independence broadcast number}, respectively.

Since the characteristic function of an independent set is a bnd-, bnr-, bn-
and h-independent broadcast, it follows from Definitions \ref{bn-i} --
\ref{h-i} that%
\begin{equation}
\alpha(G)\leq\alpha_{\operatorname{bnd}}(G)\leq\alpha_{\operatorname{bnr}%
}(G)\leq\alpha_{\operatorname{bn}}(G)\leq\alpha_{h}(G) \label{sequence2}%
\end{equation}
for any graph $G$.

When two parameters $\pi$ and $\pi^{\prime}$ are incomparable, we denote this
fact by $\pi\diamond\pi^{\prime}$. For the path $P_{n}$, where $n\geq4$, it is
easy to see that $\Gamma_{b}(P_{n})=\operatorname{IR}_{b}(P_{n}%
)=\operatorname{diam}(P_{n})=n-1$, whereas $\alpha_{h}(P_{n})=2(n-2)>\Gamma
_{b}(P_{n})$. On the other hand, for the grid graph $G_{n,n}=P_{n}\boksie
P_{n}$, if $n$ is large enough, then $\alpha_{h}(G_{n,n})=\left\lceil
\frac{n^{2}}{2}\right\rceil $ (\cite{Bouch}; see Theorem \ref{Thm-gridsB}
below), but Mynhardt and Roux \cite{MR} showed that $\Gamma_{b}(G_{n,n}%
)=\operatorname{IR}_{b}(G_{n,n})=n(n-1)>\alpha_{h}(G_{n,n})$. Therefore
$\alpha_{h}\diamond\Gamma_{b}$ and $\alpha_{h}\diamond\operatorname{IR}_{b}$,
hence $\alpha_{h}$ does not fit neatly into the inequality chain
(\ref{sequence1}). Our definitions of boundary independent broadcasts were
partially motivated by the aim of finding a definition of broadcast
independence for which the associated parameters could be inserted in
(\ref{sequence1}). Neilson \cite{LindaD} showed that $\alpha
_{\operatorname{bn}}\diamond\Gamma_{b}$ and $\alpha_{\operatorname{bnr}%
}\diamond\Gamma_{b}$, but, since a bnd-independent broadcast is minimal
dominating, $\alpha_{\operatorname{bnd}}(G)\leq\Gamma_{b}(G)$ (strict
inequality is possible). Hence%
\begin{equation}
\operatorname{ir}_{b}(G)\leq i_{\operatorname{bnd}}(G)\leq\gamma_{b}%
(G)\leq\gamma(G)\leq i(G)\leq\alpha(G)\leq\alpha_{\operatorname{bnd}}%
(G)\leq\Gamma_{b}(G)\leq\operatorname{IR}_{b}(G) \label{sequence4}%
\end{equation}
for any graph $G$. Therefore, with bnd-independent broadcasts we have achieved
this goal.

\begin{figure}[ptb]%
\centering
\includegraphics[
height=6.0148in,
width=3.9375in
]%
{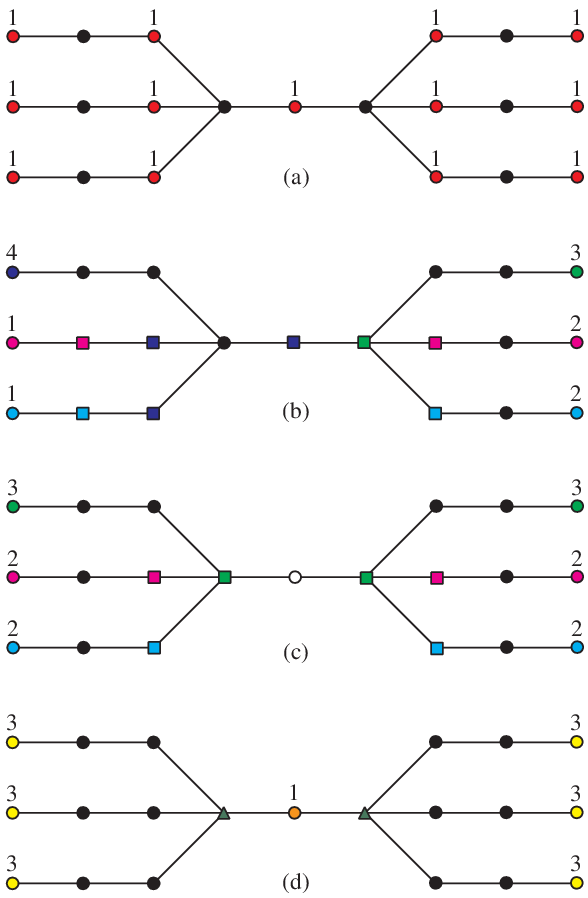}%
\caption{A tree $T$ with $\alpha(T)=\alpha_{\operatorname{bnd}}(T)<\alpha
_{\operatorname{bnr}}(T)<\alpha_{\operatorname{bn}}(T)$. A maximum independent
set is shown in (a), maximum bnd-broadcasts with weight $13$ in (a) and (b), a
maximum bnr-broadcast with weight $14$ in (c), and a maximum bn-broadcast with
weight $19$ in (d). In (b) and (c), vertices in private boundaries of
broadcasting vertices are indicated by squares, and in (d), vertices in shared
boundaries by triangles.}%
\label{fig_different}%
\end{figure}

The graph $G$ in Figure \ref{fig_different} is an example of a tree $T$ for
which $\alpha_{\operatorname{bnd}}(T)<\alpha_{\operatorname{bnr}}%
(T)<\alpha_{\operatorname{bn}}(T)$; details can be found in \cite{LindaD}.
Broadcasting from each leaf with a strength of $5$ we obtain an h-independent
broadcast with a weight of $30$, hence $\alpha_{h}(T)\geq30>\alpha
_{\operatorname{bn}}(T)$.

For the lower parameters $i_{\operatorname{bn}}$ etc., the characteristic
function of a \textbf{maximal} independent set is not necessarily a
\textbf{maximal} bn- or h-independent broadcast. For example, consider the
path $P_{6}:v_{1},...,v_{6}$, having maximal independent set $\{v_{2},v_{5}%
\}$. This set has characteristic function $f$, where $f(v_{2})=f(v_{5})=1$ and
$f(x)=0$ otherwise. The broadcast $g=(f-\{(v_{2},1)\})\cup\{(v_{2},2)\}$ is
bn- and h-independent and it is not difficult to verify that
$i_{\operatorname{bn}}(P_{6})=i_{h}(P_{6})=3>i(P_{6})=2$. On the other hand,
the corona $K_{n}\circ K_{1}$ for any complete graph $K_{n},\ n\geq4$,
satisfies $i(K_{n}\circ K_{1})=n\geq4$ but $i_{h}(K_{n}\circ K_{1}%
),\ i_{\operatorname{bn}}(K_{n}\circ K_{1})\leq3$. Therefore $i_{h}\diamond i$
and $i_{\operatorname{bn}}\diamond i$.

Dunbar et al. \cite{Dunbar} showed that every graph has a minimum weight
dominating broadcast $f$ such that $N_{f}(u)\cap N_{f}(v)=\varnothing$ for all
$u,v\in V_{f}^{+}$. Such a broadcast is maximal bnr-independent. Since any
bnr-independent broadcast is irredundant by definition, it follows that%
\begin{equation}
\operatorname{ir}_{b}(G)\leq i_{\operatorname{bnd}}(G)\leq
i_{\operatorname{bnr}}(G)\leq\gamma_{b}(G)\leq\gamma(G)\leq i(G)
\label{sequence3}%
\end{equation}
for any graph $G$. Further, although any maximal bn-independent broadcast is
dominating (see Observation \ref{rem-max-bn} below), it is not necessarily
minimal dominating, hence it is possible that $i_{\operatorname{bn}}%
>\gamma_{b}$. Neilson \cite{LindaD} showed that $i_{\operatorname{bn}}%
(G)\leq\left\lceil \frac{4}{3}\gamma_{b}(G)\right\rceil $ for all graphs $G$.

We show in Section \ref{Sec_b-i} that $\alpha_{\operatorname{bn}}(G)\leq n-1$
for all graphs $G$ of order $n$ and characterize graphs for which equality
holds. In Section \ref{bn-indep_v_h-indep} we compare $\alpha
_{\operatorname{bn}}$ and $\alpha_{\operatorname{bnr}}$ to $\alpha_{h}$ and
prove that although the differences $\alpha_{h}-\alpha_{\operatorname{bn}}$
and $\alpha_{h}-\alpha_{\operatorname{bnr}}$ can be arbitrary, the ratios
$\alpha_{h}/\alpha_{\operatorname{bn}}$ and $\alpha_{h}/\alpha
_{\operatorname{bnr}}$ are bounded by $2$ and $3$, respectively, and that
these ratios are asymptotically best possible. We deduce that $\alpha
_{h}(G)\leq2n-5$ whenever $G$ is a connected $n$-vertex graph that is not a
path. In Section \ref{Sec_bipartite} we show that $\alpha_{\operatorname{bn}%
}(G)=\alpha_{\operatorname{bnr}}(G)=\alpha_{\operatorname{bnd}}(G)=\alpha(G)$
for any $2$-connected bipartite graph $G$.

\section{Boundary independence}

\label{Sec_b-i}Suppose $f$ is a bn-independent broadcast on a graph $G$ such
that $U_{f}\neq\varnothing$; say $u\in U_{f}$. Consider the broadcast
$g_{u}=(f-\{(u,0)\})\cup\{(u,1)\}$ and notice that if any vertex $x$ of $G$
hears $u$ as well as another vertex $v\in V_{f}^{+}$, then $x\in B_{g_{u}%
}(u)\cap B_{g_{u}}(v)$. Therefore $g_{u}$ is bn-independent and $\sigma
(g_{u})>\sigma(f)$, from which we deduce that $f$ is not maximal
bn-independent. When $U_{g_{u}}\neq\varnothing$ we can repeat this process
until we obtain a maximal bn-independent broadcast $g$, i.e., one having
$U_{g}=\varnothing$. We state this fact as an observation for referencing.

\begin{observation}
\label{rem-max-bn}Any maximal bn-independent broadcast is dominating.
\end{observation}

We use Observation \ref{rem-max-bn} to prove a necessary and sufficient
condition for a bn-independent broadcast to be maximal bn-independent.

\begin{proposition}
\label{prop-max-bn}A bn-independent broadcast $f$ on a graph $G$ is maximal
bn-independent if and only if it is dominating, and either $V_{f}^{+}=\{v\}$
or $B_{f}(v)-\operatorname{PB}_{f}(v)\neq\varnothing$ for each $v\in V_{f}%
^{+}$.
\end{proposition}

\noindent\textbf{Proof.\hspace{0.1in}}Consider a maximal bn-independent
broadcast $f$ of $G$. By Observation \ref{rem-max-bn}, $f$ is dominating.
Suppose $|V_{f}^{+}|\geq2$ and there exists a vertex $v\in V_{f}^{+}$ such
that $B_{f}(v)-\operatorname{PB}_{f}(v)=\varnothing$. Since $f(v)\leq e(v)$,
the boundary $B_{f}(v)\neq\varnothing$. Since $|V_{f}^{+}|\geq2$, there exists
a vertex $w\in V_{f}^{+}-\{v\}$. By the definition of bn-independence,
$d(v,w)>f(v)$; this implies that $f(v)<e(v)$. Hence we may increase the
strength of the broadcast from $v$ to obtain the broadcast $f^{\prime
}=(f-\{(v,f(v))\})\cup\{(v,f(v)+1)\}$. Since $B_{f}(v)-\operatorname{PB}%
_{f}(v)=\varnothing$, $B_{f}(v)\subseteq\operatorname{PB}_{f}(v)$. Hence no
vertex hears $f$ from $v$ as well as from another vertex in $V_{f}^{+}$. Thus
$f^{\prime}$ is a bn-independent broadcast such that $f^{\prime}>f$. This
contradicts the maximality of $f$. Hence, if $|V_{f}^{+}|\geq2$, then
$B_{f}(v)-\operatorname{PB}_{f}(v)\neq\varnothing$ for each $v\in V_{f}^{+}$.

Conversely, suppose $f$ is a dominating bn-independent broadcast such that
either $V_{f}^{+}=\{v\}$ or $B_{f}(v)-\operatorname{PB}_{f}(v)\neq\varnothing$
for each $v\in V_{f}^{+}$. If $V_{f}^{+}=\{v\}$, then, since $f$ is
dominating, $f(v)=e(v)$ and $f$ is maximal bn-independent by definition. Hence
assume $|V_{f}^{+}|\geq2$ and $B_{f}(v)-\operatorname{PB}_{f}(v)\neq
\varnothing$ for each $v\in V_{f}^{+}$. Consider any $v\in V(G)$ and define
$f^{\prime}=(f-\{(v,f(v))\})\cup\{(v,f(v)+1)\}$. If $v\in V_{f}^{+}$, then
$B_{f}(v)-\operatorname{PB}_{f}(v)\neq\varnothing$. Let $u\in B_{f}%
(v)-\operatorname{PB}_{f}(v)$ and let $w\in V_{f}^{+}-\{v\}$ be a vertex such
that $u\in N_{f}(w)$. Since $f$ is bn-independent, $u\in B_{f}(w)$. Then
$u\in(N_{f^{\prime}}(v)\cap N_{f^{\prime}}(w))-B_{f^{\prime}}(v)$, hence
$f^{\prime}$ is not bn-independent. If $f(v)=0$, then $v\in N_{f}(w)$ for some
$w\in V_{f}^{+}$. Then $v\notin B_{f^{\prime}}(v)$ but $v\in N_{f^{\prime}%
}(v)\cap N_{f^{\prime}}(w)$. This implies that $f^{\prime}$ is not
bn-independent.~$\blacksquare$

\subsection{Bounds on boundary independence}

In this subsection we find an upper bound on the weight of a bn-independent
broadcast on a graph $G$ in terms of the size of $G$ and the sum of the
degrees of the broadcast vertices. When $G$ is a tree, this bound immediately
gives an upper bound on $\alpha_{\operatorname{bn}}(G)$. Suppose $f$ is a bn-
or bnr-independent broadcast on $G$ and an edge $xy$ of $G$ is covered by
vertices $u,v\in V_{f}^{+}$. By the definition of covered, $\{x,y\}\nsubseteq
B_{f}(u)$ and $\{x,y\}\subseteq N_{f}(u)\cap N_{f}(v)$. This violates the
bn-independence of $f$. Hence we have the following observation.

\begin{observation}
\label{edge-covered}If $f$ is a bn- or bnr-independent broadcast on a graph
$G$, then each edge of $G$ is covered by at most one vertex in $V_{f}^{+}$.
\end{observation}

\begin{proposition}
\label{Prop_UB_for_weight} Given a connected graph $G$ of size $m$, if $f$ is
a bn-independent broadcast on $G$, then $\sigma(f)\leq m-\sum_{v\in V_{f}^{+}%
}\deg(v)+|V_{f}^{+}|$.
\end{proposition}

\noindent\textbf{Proof.\hspace{0.1in}}By Observation \ref{edge-covered}, every
edge of $G$ is covered by at most one broadcast vertex. Since $f(v)\leq e(v)$
for each $v\in V_{f}^{+}$, there is at least one vertex $x$ at distance $f(v)$
from $v$. The $f(v)$ edges along the $v-x$ geodesic are all covered by $v$, as
are the remaining $\deg(v)-1$ edges incident with $v$. Therefore each
broadcast vertex $v$ covers at least $f(v)+\deg(v)-1$ edges. Counting edges we
obtain%
\[
\sum_{v\in V_{f}^{+}}(f(v)+\deg(v)-1)\leq m,
\]
which simplifies to
\[
\sigma(f)\leq m-\sum_{v\in V_{f}^{+}}\deg(v)+|V_{f}^{+}|.~\blacksquare
\]

For a broadcast $f$ on a nontrivial tree of order $n$, $\sum_{v\in V_{f}^{+}%
}\deg(v)\geq|V_{f}^{+}|$, hence the bound in Proposition
\ref{Prop_UB_for_weight} simplifies to the following bound for trees.

\begin{corollary}
\label{Cor_UB_tree}If $T$ is a tree of order $n\geq2$, then $\alpha
_{\operatorname{bnd}}(T)\leq\alpha_{\operatorname{bnr}}(T)\leq\alpha
_{\operatorname{bn}}(T)\leq n-1$.
\end{corollary}

Broadcasting from a single leaf to the whole path, it is easy to see
that$\ \alpha_{\operatorname{bnd}}(P_{n})=\alpha_{\operatorname{bnr}}%
(P_{n})=\alpha_{\operatorname{bn}}(P_{n})=n-1$ for any path $P_{n}$.

Let $f$ be an $\alpha_{\operatorname{bn}}$-broadcast on a graph $G$ and let
$T$ be a spanning tree of $G$. Removing the edges in $E(G)-E(T)$ does not
affect bn-independence, hence $f$ is also a bn-independent broadcast on $T$.
Therefore $\alpha_{\operatorname{bn}}(T)\geq\alpha_{\operatorname{bn}}(G)$,
and the result below follows from Corollary \ref{Cor_UB_tree}.

\begin{corollary}
\label{Cor_UB_graph}For any connected graph $G$ of order $n\geq2$,%
\[
\alpha_{\operatorname{bn}}(G)\leq\min\{\alpha_{\operatorname{bn}}(T):T\text{
is a spanning tree of }G\}\leq n-1.
\]

\end{corollary}

The proof of Proposition \ref{Prop_UB_for_weight} also shows that
$\sigma(f)=n-1$ if and only if every vertex in $V_{f}^{+}$ is a leaf and the
edge sets of the subtrees induced by the $f$-neighbourhoods form a partition
of $E(T)$. We use this observation to characterize graphs of order $n$ for
which $\alpha_{\operatorname{bn}}=n-1$. This characterization involves a class
of trees called spiders. As we also use spiders to show in Section
\ref{Sec_diff} that the differences $\alpha_{h}-\alpha_{\operatorname{bn}}$,
$\alpha_{h}-\alpha_{\operatorname{bnr}}$ and $\alpha_{\operatorname{bn}%
}-\alpha_{\operatorname{bnr}}$ can be arbitrary, in Section \ref{Sec_Ratio}
that the bounds for the ratios $\alpha_{\operatorname{bn}}/\alpha
_{\operatorname{bnr}},\ \alpha_{h}/\alpha_{\operatorname{bn}}$ and $\alpha
_{h}/\alpha_{\operatorname{bnd}}$ are asymptotically best possible, and in
Section \ref{Sec_Bounds} to prove a bound for $\alpha_{h}$, we define these
graphs and present results on their broadcast independence numbers in the next subsection.

\subsection{Spiders}

\label{Sec_Spiders}For $k\geq3$ and $n_{i}\geq1,\ i\in\{1,...,k\}$, the
\emph{(generalized) spider} $\operatorname{Sp}(n_{1},...,n_{k})$ is the tree
which has exactly one vertex $b$, called the \emph{head}, having $\deg(b)=k$,
and for which the $k$ components of $\operatorname{Sp}(n_{1},...,n_{k})-b$ are
paths of lengths $n_{1}-1,...,n_{k}-1$, respectively. The \emph{legs }%
$L_{1},...,L_{k}$ of the spider are the paths from $b$ to the leaves. Let
$t_{i}$ be the leaf of $L_{i},\ i=1,...,k$. If $n_{i}=r$ for each $i$, we
write $\operatorname{Sp}(r^{k})$ for $\operatorname{Sp}(n_{1},...,n_{k})$.

\begin{corollary}
\label{Cor=tree}If $G$ is a connected graph of order $n\geq2$, then
$\alpha_{\operatorname{bn}}(G)=n-1$ if and only if $G$ is a path or a spider.
\end{corollary}

\noindent\textbf{Proof.\hspace{0.1in}}Let $f$ be a bn-independent broadcast on
$G$ and assume first that $G$ is a tree. As shown in the proof of Proposition
\ref{Prop_UB_for_weight}, $\sigma(f)=n-1$ if and only if all edges of $G$ are
covered by $f$ and the number of edges covered by $v$ equals $f(v)$ for each
$v\in V_{f}^{+}$. This holds if and only if

\begin{enumerate}
\item[(1)] each $v\in V_{f}^{+}$ is a leaf and the subgraph induced by
$N_{f}(v)$ is a path of length $f(v)$.
\end{enumerate}

\noindent Since $G$ is connected and $f$ is bn-independent,

\begin{enumerate}
\item[(2)] the subpaths induced by $N_{f}(v)$ for each $v\in V_{f}^{+}$ all
have exactly one vertex in common, namely their non-broadcasting leaf.
\end{enumerate}

\noindent This is possible if and only if $G$ is a path or a generalized spider.

Now assume that $G$ has a cycle and that $\alpha_{\operatorname{bn}}(G)=n-1$.
If $G$ has a spanning tree which is not a Hamiltonian path or a spider, then
the above result for trees and Corollary \ref{Cor_UB_graph} imply that
$\alpha_{\operatorname{bn}}(G)<n-1$, which is not the case. Suppose $G$ has a
Hamiltonian path $P:v_{1},...,v_{n}$. Since $G$ has a cycle, $v_{i}v_{j}\in
E(G)$ for some $i,j$ such that $j\geq i+2$. Now $T=(P-v_{i}v_{i+1})+v_{i}%
v_{j}$ is a spanning tree of $G$ that is not a path. Since $\alpha
_{\operatorname{bn}}(G)=n-1$, we may assume that $T$ is a spider, otherwise we
have a contradiction as above.

Assume therefore that $G$ has a spanning spider $S=\operatorname{Sp}%
(n_{1},...,n_{k})$ (with notation as defined above). Consider any
$\alpha_{\operatorname{bn}}(G)$-broadcast $f$ on $G$ and let $f^{\prime}$ be
the restriction of $f$ to $S$. Then $\sigma(f^{\prime})=\sigma(f)=n-1$ and by
(1) and (2), $V_{f^{\prime}}^{+}=V_{f}^{+}=\{t_{1},...,t_{k}\}$ and
$f(t_{i})=n_{i}$ for each $i$. Since $G$ has a cycle, there is an edge $uw\in
E(G)-E(S)$. If $u$ and $w$ belong to the same leg $L_{i}$ of $S$, then
$d_{G}(t_{i},b)<f(t_{i})$, thus edges of $L_{j},\ j\neq i$, hear $f$ from both
$t_{i}$ and $t_{j}$. If $u$ and $w$ belong to different legs $L_{i},L_{j}$,
then $uw$ hears $f$ from both $t_{i}$ and $t_{j}$. Both instances contradict
$f$ being bn-independent.

We deduce that if $G$ is not a tree, then $\alpha_{\operatorname{bn}}(G)\leq
n-2$.~$\blacksquare$

\bigskip

It follows from a result in \cite{Dunbar} that $\alpha_{h}(\operatorname{Sp}%
(r^{k}))=k(2r-1)$. By Corollary \ref{Cor=tree}, $\alpha_{\operatorname{bn}%
}(\operatorname{Sp}(r^{k}))=kr$, and Neilson \cite[special case of Proposition
2.3.8]{LindaD} showed that $\alpha_{\operatorname{bnr}}(\operatorname{Sp}%
(r^{k}))\geq\alpha_{\operatorname{bnd}}(\operatorname{Sp}(r^{k}))\geq kr-k+1$.
Although there are spiders, for example $\operatorname{Sp}(1,n_{2},n_{3})$,
where $n_{2},n_{3}\geq2$, whose bnd- and bnr-independence numbers exceed
Neilson's general lower bound, it follows from our next proposition that
$\alpha_{\operatorname{bnr}}(\operatorname{Sp}(r^{k}))=\alpha
_{\operatorname{bnd}}(\operatorname{Sp}(r^{k}))=kr-k+1$ when $r\geq2$ and
$k\geq3$.

\begin{proposition}
\label{Prop_alpha_bnr_spider} If $S=\operatorname{Sp}(n_{1},...,n_{k})$ is a
spider of order $n=\sum_{i=1}^{k}n_{i}+1$, where $k\geq3$ and $n_{i}\geq2$ for
each $1\leq i\leq k$, then $\alpha_{\operatorname{bnd}}(S)=\alpha
_{\operatorname{bnr}}(\operatorname{S})=n-k$.
\end{proposition}

\noindent\textbf{Proof.\hspace{0.1in}}Again we follow the notation for spiders
as defined above. Define a broadcast $g$ on $S$ by $g(t_{1})=n_{1}$,
$g(t_{i})=n_{i}-1$ for $2\leq i\leq k$, and $g(x)=0$ otherwise. Notice that
$g$ is a dominating broadcast and $\sigma(g)=n-k$. No broadcasting vertex of
$g$ overdominates $b$ and there is exactly one broadcasting vertex on each
path $L_{i}$ for $1\leq i\leq k$. Hence $g$ is bn-independent. Further,
$\operatorname{PB}_{g}(t_{1})=\{b\}$ and for $2\leq i\leq k$, the private
boundary of $t_{i}$ consists of the vertex adjacent to $b$ on the path $L_{i}%
$. Hence $g$ is a bnr-independent and dominating broadcast. It follows that
$\alpha_{\operatorname{bnr}}(S)\geq\alpha_{\operatorname{bnd}}(S)\geq n-k$.

For the opposite inequality, let $\mathcal{F}$ be the set of $\alpha
_{\operatorname{bnr}}$-broadcasts on $S$ that minimize the number of non-leaf
broadcasting vertices. We claim that there exists a broadcast in $\mathcal{F}$
such that $b$ is not overdominated. Suppose this is not the case and consider
any $f\in\mathcal{F}$. Since $f$ is bn-independent and $b$ is overdominated,
$b$ hears $f$ from exactly one vertex $v\in V_{f}^{+}$, where possibly $v=b$.
Since $f(v)\leq e(v)$, $B_{f}(v)\neq\varnothing$. We consider two cases,
depending on whether there exists a vertex $v^{\prime}\in B_{f}(v)$ such that
$v$ and $v^{\prime}$ belong to the same leg of $S$ or not.\smallskip

\noindent\textbf{Case 1:}\hspace{0.1in}there exists a vertex $v^{\prime}\in
B_{f}(v)$ such that $v$ and $v^{\prime}$ belong to the same leg of $S$; say
$v,v^{\prime}\in V(L_{1})$. (This includes the case where $v=b$, as $b$
belongs to each leg.) Since $v$ overdominates $b$, $d(v,b)<d(v,v^{\prime})$
and $v\neq t_{1}$. Say $V(L_{1})\cap V_{f}^{+}=\{v,u_{1},...,u_{\ell}\}$.
Define the broadcast $f_{1}$ by $f_{1}(t_{1})=2f(v)-1+\sum_{i=1}^{\ell}%
f(u_{i})$, $f_{1}(x)=0$ if $x\in V(L_{1})-\{t_{1}\}$, and $f_{1}(x)=f(x)$
otherwise. Then $N_{f_{1}}(t_{1})\cap V(L_{i})\subsetneqq N_{f}(v)\cap
V(L_{i})$ for each $i$ such that $2\leq i\leq k$, which implies that $f_{1}$
is bn-independent and $\operatorname{PN}_{f_{1}}(x)\neq\varnothing$ for each
$x\in V_{f_{1}}^{+}$; that is, $f_{1}$ is bnr-independent. However, if
$f(v)>1$, then $\sigma(f_{1})>\sigma(f)$, which is impossible. We deduce that
$f(v)=1$. Since $b$ is overdominated, the only possibility is that $v=b$. But
now $f_{1}$ is an $\alpha_{\operatorname{bnr}}$-broadcast containing fewer
non-leaf broadcasting vertices than $f$, a contradiction.\smallskip

\noindent\textbf{Case 2:}\hspace{0.1in}no vertex in $B_{f}(v)$ belongs to the
same leg as $v$; assume without loss of generality that $v\in V(L_{1})-\{b\}$
and $v^{\prime}\in V(L_{2})-\{b\}$. Observe that $f(v)\geq2$ and $v$
overdominates $t_{1}$. This implies that $V(L_{1})\cap V_{f}^{+}=\{v\}$ and
also that some vertex of $L_{i}$, where $i>1$, belongs to $\operatorname{PB}%
_{f}(v)$. We may assume that $v^{\prime}\in\operatorname{PB}_{f}(v)$. Say
$f(v)=d(b,v)+q$, where $q>0$. For $i=2,3$, let $w_{i}$ be the vertex on
$L_{i}$ adjacent to $b$.

\begin{itemize}
\item Suppose first that $v^{\prime}\neq t_{2}$. Then the edge $e$ incident
with $v^{\prime}$ on the $v^{\prime}-t_{2}$-path is uncovered. Let
$V(L_{2})\cap V_{f}^{+}=\{u_{1},...,u_{\ell}\}$ and define $f_{2}$ by
$f_{2}(v)=f(v)-q$, $f_{2}(t_{2})=\sum_{i=1}^{\ell}f(u_{i})+q$, $f_{1}(x)=0$ if
$x\in V(L_{2})-\{t_{2}\}$, and $f_{2}(x)=f(x)$ otherwise. Note that
$\sigma(f_{2})=\sigma(f)$. Since $e$ is uncovered, $b\in\operatorname{PB}%
_{f_{2}}(v)$ and some vertex on the $w_{2}-t_{2}$ path belongs to
$\operatorname{PB}_{f_{2}}(t_{2})$; furthermore, $\operatorname{PB}_{f_{2}%
}(x)\supseteq\operatorname{PB}_{f}(x)$ for all $x\in V_{f_{2}}^{+}-(\{v\}\cup
V(L_{2}))$. It follows that $f_{2}$ is an $\alpha_{_{\operatorname{bnr}}}%
$-broadcast such that $b$ is not overdominated, contrary to our assumption.

\item Now suppose that $v^{\prime}=t_{2}$. Then $q=n_{2}$. Since $n_{2}\geq2$,
$q\geq2$. Let $V(L_{3})\cap V_{f}^{+}=\{u_{1},...,u_{\ell}\}$ and define
$f_{3}$ by $f_{3}(v)=f(v)-q$, $f_{3}(t_{2})=q-1$, $f_{3}(t_{3})=\sum
_{i=1}^{\ell}f(u_{i})+1$, $f_{1}(x)=0$ if $x\in V(L_{3})-\{t_{3}\}$, and
$f_{3}(x)=f(x)$ otherwise. As for $f_{2}$, $\sigma(f_{3})=\sigma(f)$. Clearly,
$b\notin N_{f_{3}}(t_{2})$, and since $q\geq2$, $b\notin N_{f_{3}}(t_{3})$.
Therefore $b\in\operatorname{PB}_{f_{3}}(v)$, $w_{2}\in\operatorname{PB}%
_{f_{3}}(t_{2})$, some vertex on the $w_{3}-t_{3}$ path belongs to
$\operatorname{PB}_{f_{3}}(t_{3})$, and $\operatorname{PB}_{f_{3}}%
(x)\supseteq\operatorname{PB}_{f}(x)$ for all $x\in V_{f_{3}}^{+}-(\{v\}\cup
V(L_{2})\cup V(L_{3}))$. As in the case of $f_{2}$, it follows that $f_{3}$ is
an $\alpha_{_{\operatorname{bnr}}}$-broadcast such that $b$ is not
overdominated, contrary to our assumption.
\end{itemize}

This completes the proof of the claim. Thus, let $f$ be an $\alpha
_{\operatorname{bnr}}$-broadcast on $S$ such that $b$ is not overdominated.
(Possibly, $b$ is not dominated at all.) Then $f(b)=0$. If $b$ is
$f$-dominated, we may assume without loss of generality that $b$ is dominated
by a vertex $v\in V(L_{1})\cap V_{f}^{+}$. Let $L_{1}^{\prime}=L_{1}$ and
$L_{i}^{\prime}=L_{i}-\{b\}$ for each $2\leq i\leq k$. Note that these paths
from a partition of $V(S)$. Restricting $f$ to each $L_{i}^{\prime}$, we
obtain $k$ separate broadcasts $f_{i}=f\upharpoonleft L_{i}^{\prime}$ for
$1\leq i\leq k$. Since $f$ is bn-independent, each $f_{i}$ is bn-independent.
Since $f$ is bnr-independent, $\operatorname{PB}_{f}(w)\neq\varnothing$ for
each $w\in V_{f}^{+}$. Also, since $b$ is not overdominated, and by the
definition of $L_{1}^{\prime}$, if $w\in V_{f}^{+}\cap V(L_{i})$, then
$\operatorname{PB}_{f}(w)\subseteq V(L_{i}^{\prime})$. Thus $\varnothing
\neq\operatorname{PB}_{f_{i}}(w)\subseteq V(L_{i}^{\prime})$. Hence the
broadcasts $f_{i}$ are bnr-independent. Since $\alpha_{\operatorname{bnr}%
}(P)=|V(P)|-1$ for any path $P$,
\[
\alpha_{\operatorname{bnd}}(S)\leq\alpha_{\operatorname{bnr}}(S)=\sigma
(f)=\sigma(f_{1})+\sum_{i=2}^{k}\sigma(f_{i})\leq n_{1}+\sum_{i=2}^{k}%
(n_{i}-1)=n-k.\ \ \blacksquare
\]

\medskip

We next determine an upper bound for $\alpha_{h}(\operatorname{Sp}%
(n_{1},...,n_{k}))$. This result generalizes the upper bound for $\alpha
_{h}(\operatorname{Sp}(r^{k}))$ in \cite{Dunbar}.

\begin{proposition}
\label{Prop_spider}If $S$ is a spider $\operatorname{Sp}(n_{1},...,n_{k})$ of
order $n$, where $k\geq3$, then $\alpha_{h}(S)\leq2n-2-k$.
\end{proposition}

\noindent\textbf{Proof.\hspace{0.1in}}Assume that $n_{1}\leq\cdots\leq n_{k}$
and note that $n=1+\sum_{i=1}^{k}n_{i}$. Let $f$ be an $\alpha_{h}$-broadcast
on $S$. If $|V_{f}^{+}|=1$, then $\sigma(f)\leq\operatorname{diam}(S)\leq
n-k+1<2n-2-k$ since $n>3$. Hence assume $|V_{f}^{+}|\geq2$. If the leg $L_{i}$
contains a broadcast vertex other than its leaf $t_{i}$, let $v$ be the
broadcast vertex on $L_{i}$ nearest to $t_{i}$. Then
\[
f^{\prime}=(f-\{(v,f(v)),(t_{i},f(t_{i}))\})\cup\{(v,0),(t_{i},f(t_{i}%
)+f(v)+1)\}
\]
is an h-independent broadcast such that $\sigma(f^{\prime})>\sigma(f)$, which
is impossible. Therefore $V_{f}^{+}\subseteq\{t_{1},...,t_{k}\}$. If the
leaves $t_{i}$ and $t_{j}$ are broadcasting vertices, then $\max
\{f(t_{i}),f(t_{j})\}\leq d(t_{i},t_{j})-1=n_{i}+n_{j}-1$. Let $l$ be the
smallest index such that $t_{l}\in V_{f}^{+}$. Since $|V_{f}^{+}|\geq2$, there
exists an index $l^{\prime}>l$ such that $t_{l^{\prime}}\in V_{f}^{+}$. Since
$f$ is h-independent, $t_{l^{\prime}}$ does not hear the broadcast from
$t_{l}$, so $t_{k}$ also does not hear the broadcast from $t_{l}$. This means
that $f(t_{l})\leq n_{l}+n_{k}-1$. Moreover, $f(t_{i})\leq n_{l}+n_{i}-1$ for
$i>l$. Hence
\[
\sigma(f)=\sum_{i=l}^{k}f(t_{i})=f(t_{l})+\sum_{i=l+1}^{k}f(t_{i})\leq
(n_{l}+n_{k}-1)+\sum_{i=l+1}^{k}(n_{l}+n_{i}-1).
\]
This inequality simplifies to%
\begin{equation}
\sigma(f)\leq n_{k}+n_{l}(k-l)+\sum_{i=l}^{k}n_{i}-(k-l)-1. \label{eq_bound}%
\end{equation}
If $l=1$, then, noting that $n_{1}\leq n_{i}$, (\ref{eq_bound}) becomes
$\sigma(f)\leq2\sum_{i=1}^{k}n_{i}-k=2n-2-k$. If $l>1$, then, noting also that
$n_{i}\geq1$, (\ref{eq_bound}) becomes
\begin{align*}
\sigma(f)  &  \leq2\sum_{i=l}^{k}n_{i}-(k-l)-1=2\sum_{i=1}^{k}n_{i}%
-2\sum_{i=1}^{l-1}n_{i}-(k-l)-1\\
&  \leq2\sum_{i=1}^{k}n_{i}-2(l-1)-(k-l)-1=2\sum_{i=1}^{k}n_{i}-(k+l)+1\\
&  <2\sum_{i=1}^{k}n_{i}-k=2n-2-k.
\end{align*}
Hence $\alpha_{h}(S)=\sigma(f)\leq2n-2-k$ and our proof is
complete.~$\blacksquare$

\section{Comparing $\alpha_{\operatorname{bn}}$ and $\alpha
_{\operatorname{bnr}}$ to $\alpha_{h}$}

\label{bn-indep_v_h-indep}In this section we show that the differences
$\alpha_{h}-\alpha_{\operatorname{bn}}$, $\alpha_{h}-\alpha
_{\operatorname{bnr}}$ and $\alpha_{\operatorname{bn}}-\alpha
_{\operatorname{bnr}}$ can be arbitrary, whereas the ratios $\alpha
_{\operatorname{bn}}/\alpha_{\operatorname{bnr}}$, $\alpha_{h}/\alpha
_{\operatorname{bn}}$ and $\alpha_{h}/\alpha_{\operatorname{bnr}}$ are bounded.

\subsection{The differences}

\label{Sec_diff}When $r\geq2$ and $k\geq3$, it follows from Proposition
\ref{Prop_alpha_bnr_spider} that
\begin{align*}
\alpha_{h}(\operatorname{Sp}(r^{k}))-\alpha_{\operatorname{bn}}%
(\operatorname{Sp}(r^{k}))  &  =k(2r-1)-kr=k\left(  r-1\right)  ,\\
\alpha_{h}(\operatorname{Sp}(r^{k}))-\alpha_{\operatorname{bnr}}%
(\operatorname{Sp}(r^{k}))  &  =k(2r-1)-(kr-k+1)=kr-1\\
\text{and\ }\alpha_{\operatorname{bn}}(\operatorname{Sp}(r^{k}))-\alpha
_{\operatorname{bnr}}(\operatorname{Sp}(r^{k}))  &  =kr-(kr-k+1)=k-1.
\end{align*}

Therefore the differences $\alpha_{h}-\alpha_{\operatorname{bn}}$, $\alpha
_{h}-\alpha_{\operatorname{bnr}}$ and $\alpha_{\operatorname{bn}}%
-\alpha_{\operatorname{bnr}}$ can be arbitrary.

\subsection{The ratios}

\label{Sec_Ratio}We show next that the ratios $\alpha_{\operatorname{bn}%
}/\alpha_{\operatorname{bnr}}$, $\alpha_{h}/\alpha_{\operatorname{bn}}$ and
$\alpha_{h}/\alpha_{\operatorname{bnr}}$ are bounded. When $f$ is a
bnr-broadcast, $\operatorname{PB}_{f}(v)\neq\varnothing$ for each $v\in
V_{f}^{+}$, but when $f$ is an h- or bn-independent broadcast, it is possible
that $\operatorname{PB}_{f}(v)=\varnothing$ for some $v\in V_{f}^{+}$. For
each of these three types of broadcasts, if $f(v)=1$, then $v\in
\operatorname{PB}_{f}(v)$. Therefore we have the following observation.

\begin{observation}
\label{Ob_PB}If $f$ is an h- or a bn-independent broadcast such that
$\operatorname{PB}_{f}(v)=\varnothing$ for some $v\in V_{f}^{+}$, then $v\in
V_{f}^{++}$.
\end{observation}

\begin{theorem}
\label{Thm_bnr/bn}For any graph $G$, $\alpha_{\operatorname{bn}}%
(G)/\alpha_{\operatorname{bnr}}(G)<2$, and this bound is asymptotically best possible.
\end{theorem}

\noindent\textbf{Proof.\hspace{0.1in}}Let $f$ be an $\alpha_{\operatorname{bn}%
}$-broadcast on $G$. If $\operatorname{PB}_{f}(v)\neq\varnothing$ for each
$v\in V_{f}^{+}$, then $f$ is bnr-independent and $\alpha_{\operatorname{bn}%
}(G)=\alpha_{\operatorname{bnr}}(G)$. Hence assume $\operatorname{PB}%
_{f}(v)=\varnothing$ for some $v\in V_{f}^{+}$. Then $|V_{f}^{+}|\geq2$ and,
by Observation \ref{Ob_PB}, $v\in V_{f}^{++}$. If $V_{f}^{1}=\varnothing$,
choose an arbitrary vertex $u\in V_{f}^{++}$. Define the broadcast $g$ by
\[
g(x)=\left\{
\begin{tabular}
[c]{ll}%
$f(x)-1$ & if $V_{f}^{1}=\varnothing$ and $x\in V_{f}^{++}-\{u\}$\\
$f(x)-1$ & if $V_{f}^{1}\neq\varnothing$ and $x\in V_{f}^{++}$\\
$f(x)$ & otherwise.
\end{tabular}
\ \ \ \right.
\]
Then $\sigma(g)\geq\sigma(f)-|V_{f}^{++}|\geq\frac{1}{2}\sigma(f)$ and at
least one of the inequalities is strict. Moreover, since $f$ overlaps only in
boundaries and $g(x)<f(x)$ for each $x\in V_{f}^{++}$ (if $V_{f}^{1}%
\neq\varnothing$) or for each but one $x\in V_{f}^{++}$ (if $V_{f}^{+}%
=V_{f}^{++}$), the $g$-neighbourhoods are pairwise disjoint. Since $g(x)\leq
e(x)$ for each $x\in V_{g}^{+}$, there is at least one vertex at distance
$g(x)$ from $x$. Hence $B_{g}(x)\neq\varnothing$, and since the $g$%
-neighbourhoods are pairwise disjoint, $\operatorname{PB}_{g}\neq\varnothing$
for each $x\in V_{g}^{+}$. Therefore $g$ is bnr-independent and $\alpha
_{\operatorname{bnr}}(G)\geq\sigma(g)>\frac{1}{2}\alpha_{\operatorname{bn}%
}(G)$.

To see that the bound is asymptotically best possible, consider the spiders
$S=\operatorname{Sp}(2^{k}),\ k\geq3$. Since $\alpha_{\operatorname{bn}%
}(S)=2k$ and $\alpha_{\operatorname{bnr}}(S)=k+1$ (Proposition
\ref{Prop_alpha_bnr_spider}), the result follows.~$\blacksquare$

\bigskip

We now bound $\alpha_{h}/\alpha_{\operatorname{bn}}$ and $\alpha_{h}%
/\alpha_{\operatorname{bnr}}$. Since the proofs overlap, we state the results
as parts of the same theorem.

\begin{theorem}
\label{Thm_h-bn}For any graph $G$,\vspace{-0.1in}

\begin{enumerate}
\item[$(i)$] $\alpha_{h}(G)/\alpha_{\operatorname{bn}}(G)<2$, and

\item[$(ii)$] $\alpha_{h}(G)/\alpha_{\operatorname{bnr}}(G)<3$. \vspace
{-0.08in}
\end{enumerate}

Both bounds are asymptotically best possible.
\end{theorem}

\noindent\textbf{Proof.\hspace{0.1in}}$(i)\hspace{0.1in}$Let $f$ be an
$\alpha_{h}$-broadcast on $G$. If $f$ is bn-independent, then $\alpha
_{h}(G)=\alpha_{\operatorname{bn}}(G)$ and we are done, hence assume $v,w\in
V_{f}^{+}$ cover the same edge, say $e$. Since $f$ is h-independent, no
broadcasting vertex hears any other broadcasting vertex. In particular,
neither $v$ nor $w$ is incident with $e$. Hence $v,w\in V_{f}^{++}$. Define
the broadcast $f^{\prime}$ on $G$ by $f^{\prime}(x)=\left\lceil \frac{f(x)}%
{2}\right\rceil $ if $x\in V_{f}^{++}$ and $f^{\prime}(x)=f(x)$ otherwise.

We claim that for $v,w\in V_{f}^{++}$, if at least one of $f(v)$ and $f(w)$ is
even, then no vertex of $G$ hears $f^{\prime}$ from both $v$ and $w$, while if
$f(v)$ and $f(w)$ are both odd, then $N_{f^{\prime}}(v)\cap N_{f^{\prime}%
}(w)\subseteq B_{f^{\prime}}(v)\cap B_{f^{\prime}}(w)$. This will show that
$f^{\prime}$ is bn-independent.

Suppose there exists a vertex $u\in N_{f^{\prime}}(v)\cap N_{f^{\prime}}(w)$
for some $v,w\in V_{f}^{++}$. Then $f^{\prime}(v)\geq d(v,u)$, $f^{\prime
}(w)\geq d(w,u)$ and $d(v,w)\leq f^{\prime}(v)+f^{\prime}(w)$. If $f(v)\neq
f(w)$, say without loss of generality $f(w)<f(v)$, then%
\[
d(v,w)\leq f^{\prime}(v)+f^{\prime}(w)=\left\lceil \frac{f(v)}{2}\right\rceil
+\left\lceil \frac{f(w)}{2}\right\rceil \leq f(v).
\]
But then $w\in V_{f}^{+}$ hears $v\in V_{f}^{+}-\{w\}$, contradicting the
h-independence of $f$. If $f(v)=f(w)\equiv0\ (\operatorname{mod}\ 2)$, then%
\[
d(v,w)\leq\left\lceil \frac{f(v)}{2}\right\rceil +\left\lceil \frac{f(w)}%
{2}\right\rceil =f(v),
\]
again contradicting the h-independence of $f$. Finally, if $f(v)=f(w)\equiv
1\ (\operatorname{mod}\ 2)$, then%
\[
d(v,w)\leq\left\lceil \frac{f(v)}{2}\right\rceil +\left\lceil \frac{f(w)}%
{2}\right\rceil =f(v)+1.
\]
Since $f$ is h-independent, $d(v,w)=f(v)+1=2f^{\prime}(v)=2f^{\prime}(w)$ and
$u\in B_{f^{\prime}}(v)\cap B_{f^{\prime}}(w)$. It follows that $f^{\prime}$
is bn-independent.

If $f(v)$ is odd for at least one $v\in V_{f}^{++}$, then $\alpha
_{\operatorname{bn}}(G)\geq\sigma(f^{\prime})>\frac{1}{2}\sigma(f)$. If $f(v)$
is even for each $v\in V_{f}^{++}\neq\varnothing$, then $f^{\prime}$ is not
maximal bn-independent, for at least one $f^{\prime}(v)$ can be increased
without any edge being covered by more than one vertex, and $\alpha
_{\operatorname{bn}}(G)>\sigma(f^{\prime})\geq\frac{1}{2}\sigma(f)$. If
$V_{f}^{+}=V_{f}^{1}$, then $\alpha_{\operatorname{bn}}(G)=\alpha_{h}(G)$.
Hence $\alpha_{h}(G)/\alpha_{\operatorname{bn}}(G)<2$.

$(ii)\hspace{0.1in}$If every vertex of $G$ hears $f^{\prime}$ (as defined
above) from exactly one vertex in $V_{f^{\prime}}^{+}$, then $f^{\prime}$ is a
bnr-independent broadcast and we are done, hence assume that a vertex $u$
hears $f^{\prime}$ from two vertices $v$ and $w$. Since $f^{\prime}$ is
bn-independent, $u\in B_{f^{\prime}}(v)\cap B_{f^{\prime}}(w)$. From the
analysis above, this happens if and only if $v,w\in V_{f}^{++}$ and
$f(v)=f(w)\equiv1\ (\operatorname{mod}\ 2)$. Therefore $f(v),f(w)\geq3$.
Choose any vertex $z\in V_{f}^{++}$ such that $f(z)$ is odd. Define the
broadcast $f^{\prime\prime}$ by
\[
f^{\prime\prime}(x)=\left\{
\begin{tabular}
[c]{cl}%
$\left\lceil \frac{f(x)}{2}\right\rceil $ & if $x=z$\\
& \\
$\left\lfloor \frac{f(x)}{2}\right\rfloor $ & if $x\in V_{f}^{++}-\{z\}$\\
& \\
$f(x)$ & otherwise.
\end{tabular}
\ \ \ \ \right.
\]
Then $N_{f^{\prime\prime}}(v)\cap N_{f^{\prime\prime}}(w)=\varnothing$ for all
$v\in V_{f}^{++}$ and $w\in V_{f}^{+}$, hence $f^{\prime\prime}$ is
bnr-independent. Moreover, $\sigma(f^{\prime\prime})>\sigma(f)-\frac{2}%
{3}\sigma(f)=\frac{1}{3}\sigma(f)$. Hence $\alpha_{\operatorname{bnr}}%
(G)\geq\sigma(f^{\prime\prime})>\frac{1}{3}\sigma(f)=\frac{1}{3}\alpha_{h}%
(G)$, i.e., $\alpha_{h}(G)<3\alpha_{\operatorname{bnr}}(G)$.

The spiders $\operatorname{Sp}(r^{k})$, which satisfy
\[
\alpha_{h}(\operatorname{Sp}(r^{k}))=k(2r-1)\text{\ and\ }\alpha
_{\operatorname{bn}}(\operatorname{Sp}(r^{k}))=kr,
\]
show that the ratio $\alpha_{h}/\alpha_{\operatorname{bn}}<2$ is
asymptotically best possible. The spiders $\operatorname{Sp}(2^{k})$, which
satisfy $\alpha_{h}(\operatorname{Sp}(2^{k}))=3k$ and $\alpha
_{\operatorname{bnr}}(\operatorname{Sp}(2^{k}))=k+1$, illustrate the
corresponding result for the ratio $\alpha_{h}/\alpha_{\operatorname{bnr}}%
<3$.~$\blacksquare$

\subsection{Bounds}

\label{Sec_Bounds}Theorem \ref{Thm_h-bn} and any upper bounds for
$\alpha_{\operatorname{bn}}$ or $\alpha_{\operatorname{bnr}}$ can be used to
obtain upper bounds for $\alpha_{h}$. Conversely, lower bounds for $\alpha
_{h}$ provide lower bounds for $\alpha_{\operatorname{bn}}$ and $\alpha
_{\operatorname{bnr}}$. Bessy and Rautenbach \cite{BR} obtained a general
upper bound for $\alpha_{h}$. For a broadcast $f$ on $G$, define $f_{\max
}=\max\{f(v):v\in V(G)\}$.

\begin{theorem}
\label{Thm_BR}\emph{\cite{BR}}\hspace{0.1in}If $G$ is a connected graph such
that
\[
\max\{\operatorname{diam}(G),\alpha(G)\}\geq3,
\]
and $f$ is a maximal h-independent broadcast on $G$, then
\[
\sigma(f)\leq4\alpha(G)-4\min\left\{  1,\frac{2\alpha(G)}{f_{\max}+2}\right\}
.
\]

\end{theorem}

Therefore $\alpha_{h}(G)<4\alpha(G)$, giving the ratio $\alpha_{h}%
(G)/\alpha(G)<4$ whenever $G$ satisfies the conditions of Theorem
\ref{Thm_BR}. The bound on the ratio is asymptotically best possible, since
$\alpha_{h}(P_{n})=2(n-2)$ when $n\geq4$, whereas $\alpha(P_{n})=\left\lceil
n/2\right\rceil $.

We present a sharp upper bound for $\alpha_{h}(G)$ in terms of the order of
$G$ as a corollary to our previous results.

\begin{corollary}
\label{Cor_UB-h}If $G$ is a connected graph of order $n$ that is not a path,
then $\alpha_{h}(G)\leq2n-5$.
\end{corollary}

\noindent\textbf{Proof.\hspace{0.1in}}When $G$ is not a spider, the result
follows immediately from Corollaries \ref{Cor_UB_graph} and \ref{Cor=tree} and
Theorem \ref{Thm_h-bn}$(i)$. By Proposition \ref{Prop_spider}, $\alpha
_{h}(\operatorname{Sp}(n_{1},...,n_{k}))$\ $\leq2n-2-k\leq2n-5$ when $k\geq
3$.~$\blacksquare$

\bigskip

Since $\operatorname{Sp}(r^{3})$ has order $3r+1$ and $\alpha_{h}%
(\operatorname{Sp}(r^{3}))=3(2r-1)=2(3r+1)-5$, the bound in Corollary
\ref{Cor_UB-h} is sharp. For graphs with large independence numbers, this
bound is better than the bound in Theorem \ref{Thm_BR}. If $G\neq P_{n}$ is a
connected graph of order $n$ such that $\alpha(G)=(1-\varepsilon)n$, where
$\varepsilon\leq\frac{1}{2}$ (which is the case when $G$ is bipartite, for
example), then Corollary \ref{Cor_UB-h} gives
\[
\alpha_{h}(G)\leq2n-5=\frac{2\alpha(G)}{1-\varepsilon}-5<4\alpha
(G)-4\min\left\{  1,\frac{2\alpha(G)}{f_{\max}+2}\right\}  .
\]

Erwin \cite{Ethesis} noted that if a connected graph $G$ has order $n\geq4$,
then any $\alpha_{h}$-broadcast on $G$ has $|V_{f}^{+}|\geq2$. Broadcasting
from two antipodal vertices $v,w$ such that $f(v)=f(w)=\operatorname{diam}%
(G)-1$, Erwin therefore obtained that $\alpha_{h}(G)\geq2(\operatorname{diam}%
(G)-1)$. Dunbar et al. \cite{Dunbar} improved Erwin's bound as follows; note
that the bound is sharp for (e.g.) $\operatorname{Sp}(r^{k})$. Let $\mu(G)$
denote the cardinality of a largest set of mutually antipodal vertices in $G$.

\begin{proposition}
\label{Prop_mu}\emph{\cite{Dunbar}}\hspace{0.1in}If $G$ is a connected graph
$G$ order at least $3$, then $\alpha_{h}(G)\geq\mu(G)(\operatorname{diam}%
(G)-1)$, and this bound is sharp.
\end{proposition}

Theorem \ref{Thm_h-bn} and Proposition \ref{Prop_mu} immediately give the
following lower bounds for $\alpha_{\operatorname{bn}}$ and~$\alpha
_{\operatorname{bnr}}$.

\begin{corollary}
\label{Cor_lbs_bn}For any connected graph $G$ of order at least $3$,
\[
\alpha_{\operatorname{bn}}(G)\geq\frac{1}{2}\mu(G)(\operatorname{diam}%
(G)-1)+1
\]
and
\[
\alpha_{\operatorname{bnr}}(G)\geq\frac{1}{3}\mu(G)(\operatorname{diam}%
(G)-1)+1.
\]
Both bounds are sharp.
\end{corollary}

For the path $P_{n}$, where $n\geq3$, the bound for $\alpha_{\operatorname{bn}%
}$ is%
\[
\alpha_{\operatorname{bn}}(P_{n})\geq\operatorname{diam}(P_{n})=n-1,
\]
which gives the exact value for $\alpha_{\operatorname{bn}}(P_{n})$, and for
the spider $S=\operatorname{Sp}(2^{k})$, the bound for $\alpha
_{\operatorname{bnr}}$ is $\alpha_{\operatorname{bnr}}(S)\geq k+1$, which also
gives $\alpha_{\operatorname{bnr}}(S)$ exactly.

\section{Bipartite graphs}

\label{Sec_bipartite}It is well known that for the $m\times n$ grid graph
$G_{m,n}=P_{m}\boksie P_{n}$, $\alpha(G_{m,n})=\left\lceil \frac{mn}%
{2}\right\rceil $. Determining the domination number of grid graphs was a
major problem in domination theory until Chang's conjecture, $\gamma
(G_{m,n})=\left\lfloor \frac{(m+2)(n+2)}{5}\right\rfloor -4$ for $m,n$ such
that $16\leq m\leq n$ \cite{Chang}, was proved by Gon\c{c}alves, Pinlou, Rao
and Thomass\'{e} \cite{GPRT}. Therefore grid graphs form an important class of
graphs to consider for other domination parameters. Also, Bouchemakh and Zemir
\cite{Bouch} considered h-independence for grids, making it one of the few
classes of graphs for which any work on independent broadcasts had been done
prior to the dissertation \cite{LindaD}.

We prove a result for $2$-connected bipartite graphs from which we immediately
obtain $\alpha_{\operatorname{bnr}}(G_{m,n})$ and $\alpha_{\operatorname{bn}%
}(G_{m,n})$.

\begin{theorem}
\label{Thm_bn-bipartite}If $G$ is a $2$-connected bipartite graph, then
\[
\alpha_{\operatorname{bn}}(G)=\alpha_{\operatorname{bnr}}(G)=\alpha
_{\operatorname{bnd}}(G)=\alpha(G).
\]

\end{theorem}

\noindent\textbf{Proof.\hspace{0.1in}}We prove that $G$ has an $\alpha
_{\operatorname{bn}}$-broadcast $f$ such that $f(v)=1$ for each $v\in
V_{f}^{+}$. Among all $\alpha_{\operatorname{bn}}$-broadcasts of $G$, let $f$
be one for which $|V_{f}^{++}|$ is minimum. When $V_{f}^{++}=\varnothing$, we
are done, hence assume there exists $v\in V_{f}^{++}$. Let $f(v)=k\geq2$.
Since $f(v)\leq e(v)$, there is a vertex $u$ at distance $k$ from $v$. Since
$G$ is $2$-connected, $u$ and $v$ lie on a common cycle; let $C$ be the
shortest cycle containing $u$ and $v$. Suppose $u$ is the only vertex such
that $d(u,v)=k$. Then $C$ has length $2k$. Let $C:v=v_{0},v_{1},...,v_{2k}=v$.
Define the broadcast $g$ by
\[
g(x)=\left\{
\begin{tabular}
[c]{ll}%
$0$ & if $x=v_{i}$ and $i\equiv k\ (\operatorname{mod}\ 2)$\\
$1$ & if $x=v_{i}$ and $i\equiv k+1\ (\operatorname{mod}\ 2)$\\
$f(x)$ & otherwise.
\end{tabular}
\ \ \right.
\]
Suppose there is another vertex $w\neq u$ at distance $k$ from $v$ on $C$.
Then there is a $u-w$ path of length $2k$ containing $v$, say $P:u=v_{0}%
,v_{1},...,v_{k}=v,...,v_{2k}=w$. Define $g$ by%
\[
g(x)=\left\{
\begin{tabular}
[c]{ll}%
$0$ & if $x=v_{i}$ and $i\equiv0\ (\operatorname{mod}\ 2)$\\
$1$ & if $x=v_{i}$ and $i\equiv1\ (\operatorname{mod}\ 2)$\\
$f(x)$ & otherwise.
\end{tabular}
\ \ \right.
\]
In either case, since $G$ is bipartite, no two vertices $v_{i},v_{j}$ (on $P$
or $C$) where $i\equiv j\ (\operatorname{mod}\ 2)$ are adjacent. Also,
$N_{g}(v_{i})\subseteq N_{f}(v)$ for each $i$. Hence $g$ is bn-independent.
Notice that $\sigma(g)=\sigma(f)$. Thus either $g$ contradicts the minimality
of $|V_{f}^{++}|$ among the $\alpha_{\operatorname{bn}}$-broadcasts of $G$, or
$g$ is not maximal bn-independent and contradicts $f$ being an $\alpha
_{\operatorname{bn}}$-broadcast.

Hence $G$ has an $\alpha_{\operatorname{bn}}$-broadcast $f$ such that $f(v)=1$
for each $v\in V_{f}^{+}$. Then $V_{f}^{1}$ is an independent set, from which
we deduce that $\alpha_{\operatorname{bn}}(G)\leq\alpha(G)$. The result
follows from the inequalities (\ref{sequence2}).~$\blacksquare$

\bigskip

Since $\alpha_{\operatorname{bn}}(\operatorname{Sp}(3^{k}))=3k$,
$\alpha_{\operatorname{bnr}}(\operatorname{Sp}(3^{k}))=2k+1$ and
$\alpha(\operatorname{Sp}(3^{k}))=2k$, Theorem \ref{Thm_bn-bipartite} does not
hold for bipartite graphs that are not $2$-connected.

Bouchemakh and Zemir \cite{Bouch} determined $\alpha_{h}$ for all grid graphs,
showing that when $m$ and $n$ are large enough, $\alpha_{h}(G_{m,n}%
)=\alpha(G_{m,n})=\left\lceil \frac{mn}{2}\right\rceil $.

\begin{theorem}
\label{Thm-gridsB}\emph{\cite{Bouch}}\hspace{0.1in}$(i)\hspace{0.1in}$If
$m,n\in\mathbb{Z}$ such that $2\leq m\leq n$ and $m\leq4$, then $\alpha
_{h}(G_{m,n})=2(m+n-3)=2(\operatorname{diam}(G_{m,n})-1)$.

\begin{enumerate}
\item[$(ii)$] If $m,n\in\mathbb{Z}$ such that $5\leq m\leq n$ and
$(m,n)\notin\{(5,5),(5,6)\}$, then $\alpha_{h}(G_{m,n})=\left\lceil \frac
{mn}{2}\right\rceil .$

\item[$(iii)$] $\alpha_{h}(G_{5,5})=15$ and $\alpha_{h}(G_{5,6})=16$.
\end{enumerate}
\end{theorem}

\noindent It therefore follows from the inequalities (\ref{sequence2}) that
for $n\geq m\geq5$ and $(m,n)\notin\{(5,5),(5,6)\}$,
\[
\alpha(G_{m,n})=\alpha_{\operatorname{bnd}}(G_{m,n})=\alpha
_{\operatorname{bnr}}(G_{m,n})=\alpha_{\operatorname{bn}}(G_{m,n})=\alpha
_{h}(G_{m,n})=\left\lceil \frac{mn}{2}\right\rceil .
\]
However, Theorem \ref{Thm_bn-bipartite} immediately gives%
\[
\alpha(G_{m,n})=\alpha_{\operatorname{bnd}}(G_{m,n})=\alpha
_{\operatorname{bnr}}(G_{m,n})=\alpha_{\operatorname{bn}}(G_{m,n})=\left\lceil
\frac{mn}{2}\right\rceil
\]
whenever $m$ and $n$ are integers such that\ $2\leq m\leq n$.

\section{Future work}

\label{Sec_future}Although $i_{\operatorname{bnd}}$ and $\alpha
_{\operatorname{bnd}}$ fit nicely into the inequality chain (\ref{sequence4}),
the definition of bnd-independence forces this to be the case. The concept is
difficult to work with and not very much is known about it. For example,
although the difference $\alpha_{\operatorname{bnr}}-\alpha
_{\operatorname{bnd}}$ can be arbitrary for trees \cite{LindaD}, the behaviour
of $\alpha_{\operatorname{bnr}}/\alpha_{\operatorname{bnd}}$ has not been
determined. It would also be interesting, for comparison, to determine
$\alpha_{\operatorname{bnd}}(G)$ for classes of graphs for which $\alpha
_{h}(G)$, $\alpha_{\operatorname{bn}}(G)$ or $\alpha_{\operatorname{bnr}}(G)$
is known.

For h-independence it would be interesting to find more graphs (if they exist)
for which the bound in Corollary \ref{Cor_UB-h} is sharp.

\label{refs}

\bigskip

\bigskip

\bigskip
\end{document}